\input amssym.def 
\input amssym
\magnification=1200
\parindent0pt
\hsize=16 true cm
\baselineskip=13  pt plus .2pt
$ $

\def\Z{{\Bbb Z}}
\def\D{{\Bbb D}}
\def\A{{\Bbb A}}
\def\S{{\Bbb S}}

\centerline {\bf  Large finite group actions on surfaces: Hurwitz groups, maximal
reducible}

\centerline {\bf and maximal handlebody groups, bounding and non-bounding actions}

\bigskip

\centerline {Bruno P. Zimmermann}

\medskip

\centerline {Universit\`a degli Studi di Trieste}

\centerline {Dipartimento di Matematica e Geoscienze}

\centerline {34127 Trieste, Italy}

\bigskip \bigskip

{\bf Abstract.}  We consider large finite group-actions on surfaces and discuss and
compare various notions for such actions: Hurwitz actions and Hurwitz groups; 
maximal reducible and completely reducible actions; bounding and geometrically
bounding actions; maximal handlebody groups and maximal bounded surface groups; in
particular, we discuss small simple groups of various types. 

\medskip

A  Hurwitz group is a finite group of
orientation-preserving diffeomorphisms of maximal possible order $84(g-1)$ of
a closed orientable surface of genus $g>1$.  A maximal handlebody group
instead  is a group of orientation-preserving diffeomorphisms  of maximal
possible order $12(g-1)$  of a 3-dimensional handlebody of genus $g>1$. Among others,
we consider the question of when a Hurwitz group acting on a surface of genus $g$
contains a subgroup of maximal possible order $12(g-1)$ extending  to a handlebody
(or, more generally, a maximal reducible group extending to   a product with
handles), and show that such Hurwitz groups are closely related to the smallest
Hurwitz group ${\rm PSL}_2(7)$ of order 168 acting on Klein's quartic of genus 3. We
discuss simple groups of small order which are maximal handlebody groups and, more
generaly, maximal reducible groups. We discuss also the problem of which Hurwitz
actions bound geometrically, and in particular whether Klein's quartic bounds
geometrically: does there exist a compact hyperbolic 3-manifold with totally
geodesic boundary isometric to Klein's quartic?  Finally, large bounding and
non-bounding actions on surfaces of genus 2, 3  and 4 are discussed in section 3.

\bigskip \bigskip

{\bf 1. Introduction}

\medskip

All finite group actions in the present paper will be
orientation-preserving, all manifolds will be orientable.

\bigskip

1.1  {\it Hurwitz groups}

\medskip

By the formula
of Riemann-Hurwitz, the maximum order of a finite group $H$ of
diffeomorphisms of a closed surface
$\Sigma$ of genus $g>1$ is
$84(g-1)$; such a group $H$ is called a {\it Hurwitz groups 
of genus $g$}. The quotient orbifold $\Sigma/H$ 	of such a {\it Hurwitz
action} is the 2-sphere with three branch points of orders 2, 3 and 7, and
the Hurwitz groups are exactly the finite quotients of its orbifold
fundamental group which is isomorphic to the  hyperbolic triangle group of type
(2,3,7), with presentation 
$$<x, y, z \mid  x^2 = y^3 = z^7  = xyz = 1>.$$
This triangle group acts by isometries on the hyperbolic plane $\Bbb H^2$,
there is an exact sequence
$$1 \to K \hookrightarrow  (2,3,7)  \to H \to 1,$$
and the factor group $H \cong (2,3,7)/K$ acts by isometries on the hyperbolic
surface $\Bbb H^2/K$.

\bigskip

1.2  {\it  Reducible and irreducible actions; products with handles}

\medskip

An action of a finite group $G$ on a surface $\Sigma$ is 
{\it reducible} if there is a nontrivial simple closed curve on $\Sigma$
which is mapped either to itself or to a disjoint curve by every element of
$G$, and otherwise {\it irreducible}. For a reducible action, one can equivariantly
attach thickened 2-disks along the curve and its $G$-images to the boundary
component $\Sigma \times 1$ of the product
$\Sigma \times I$  to obtain a 3-manifold which is a
product with handles of genus $g$ to which the $G$-action extends. Iterating
the construction  one finishes with a product with handles ${\cal P}$ with a
$G$-action such that the action of each stabilizer of a boundary component of
${\cal P}$ arising from  $\Sigma \times 1$ is irreducible.

\medskip

In general, a  {\it product with handles} $\cal P$ is 
defined as follows. Considering the product $\Sigma \times I$ of a
closed surface $\Sigma$ of genus $g$ with the unit interval, one attaches
2-handles
$D^2 \times I$ (where $D^2$ denotes the 2-disk) along the boundary parts
$S^1 \times I$ to a finite collection of disjoint simple closed curves on
the surface $\Sigma \times 1$, then one closes eventually created
2-sphere boundary components by attaching 3-balls to obtain a
3-manifold $\cal P$ whose boundary consists of the surface 
$\partial_0  {\cal P} =  \Sigma \times 0$ of genus $g$ (the {\it outer boundary} of
$\cal P$) and a collection of surfaces $\partial_1 {\cal P}$ (the {\it inner
boundary} which may be empty). We call $\cal P$ a {\it product  with handles of
genus $g$} (see section 4 for an explicit example).

\medskip

A $G$-action on a surface $\Sigma$ is irreducible if and only if the quotient
orbifold $\Sigma/G$ is the 2-sphere with exactly three branch points. In
particular, all Hurwitz actions are irreducible.
By the formula of Riemann-Hurwitz, the maximal possible order of a reducible
action of $G$ on a surface of genus $g>1$ is $12(g-1)$, and in this
case the quotient orbifold $\Sigma/G$ is the 2-sphere with four branch
points of orders 2, 2, 2 and 3 whose orbifold fundamental group is isomorphic
to the hyperbolic quadrangle group of type (2,2,2,3), with presentation
$$<x_1,x_2,x_3,x_4 \mid  x_1^2 = x_2^2 = x_3^2 = x_4^3 = x_1x_2x_3x_4 = 1>.$$
As before, there is an exact sequence

$$1 \to K \hookrightarrow  (2,2,2,3) \to G \to 1,$$

and the factor group $G \cong (2,2,2,3)/K$ acts by isometries on the
hyperbolic surface
$\Bbb H^2/K$. Let $n$ denote the order of the  image of the element $x_1x_2$
in the finite group $G$; by adding the relation
$(x_1x_2)^n = 1$ to the above presentation, we obtain a presentation of the
free poduct with amalgamation

$$G_n \; = \; \D_n*_{\Z_n}(2,3,n)$$

where $(2,3,n)$ denotes the triangle group of type $(2,3,n)$ and $\D_n =
(2,2,n)$ the dihedral group of order $2n$. So the finite groups $G$ 
admitting a reducible action of maximal possible order $12(g-1)$ on a surface
of genus $g>1$ are exactly the finite factor groups with torsionfree
kernel of the quadrangle group (2,2,2,3), or equivalently  of one of the
groups $G_n$, $n \ge 2$; we call such a group a 
{\it maximal reducible group}, and specifically a
{\it $G_n$-group}. So the maximal reducible groups are
exactly the  $G_n$-groups, for some $n \ge 2$.

\medskip

A quotient of a product with handles by a finite group action is an {\it
orbifold product with handles}; such an orbifold can be uniformized by a function
group, i.e. by a Kleinian groups with an invariant component of the regular set. An
approach to orbifold product with handles and the construction of the
uniformizing Kleinian groups is given in [R], using the language of finite  graphs of
groups.

\bigskip

1.3  {\it  Completely reducible actions, maximal handlebody groups and maximal
bounded surface groups}

\medskip

An action of a finite group $G$ on a surface $\Sigma_g$ is {\it completely
reducible} if it extends to a handlebody $V_g$ with $\partial V_g =
\Sigma_g$. By the equivariant Dehn Lemma/loop theorem [MY], a completely
reducible action is reducible, and the maximum order of a completely
reducible action of a finite group $G$ acting on  a 3-dimensional
handlebody $V_g$ of genus $g>1$ is again $12(g-1)$ (as first observed
in [Z1], prior to the appearence of the equivariant Dehn lemma/loop theorem,
see also [MMZ, Theorem 7.2]); such a group $G$ is called a {\it maximal
handlebody group of genus $g$}. By [Z2] (see also [Z5], [MMZ] for the
following), the maximal handlebody groups  $G$ are  exactly the finite
quotients with torsionfree kernel of one of the four groups

$$G_2 = \D_2*_{\Z_2}\D_3, \;\;\;  G_3 = \D_3*_{\Z_3}{\A_4}, \;\;\;
G_4 = \D_4*_{\Z_4}{\S_4}, \;\;\;  G_5 = \D_5*_{\Z_5}{\A_5}$$

where $\A_4 =(2,3,3)$ denotes the tetrahedral group
of order 12,
$\S_4 = (2,3,4)$ the octahedral group of order 24 and $\A_5 = (2,3,5)$ the
dodecahedral group of order 60; together with the dihedral groups $\D_n =
(2,2,n)$ of order $2n$, these are exactly the finite or spherical triangle
groups and, together with the cyclic groups, the finite subgroups of SO(3).

\medskip

The quotient orbifold $V_g/G$ of a maximal handlebody group $G$ of genus
$g$ by a  finite group action is a {\it handlebody orbifold} whose orbifold
fundamental group is one of the four groups $G_2, G_3, G_4$ or $G_5$; the
underlying topological space is the 3-disk, its boundary
$\partial V_g/G$  is the 2-sphere with four branch points of orders 2, 2, 2
and 3.  The handlebody  orbifolds are uniformized by virtual Schottky groups
(i.e., virtually free function groups with a single component of the regular set),
cf.  [MMZ], [Z5].

\bigskip

An upper bound for the  maximal order of a finite group of
homeomorphisms of a compact surface with nonempty boundary (orientable or
not) of algebraic genus
$g$ (the rank of its free fundamental group) is also $12(g-1)$. This can be
seen by taking the product of the surface with a closed interval (twisted if
the surface is nonorientable); the result is an orientable handlebody of
genus $g$ to which the finite group action extends orientation-preservingly, so
the upper bound for handlebodies applies, and such groups of order $12(g-1)$
are called {\it maximal bounded surface groups}. The maximal bounded surface groups
are exactly the finite quotients with torsionfree kernel of the group 

$$G_2 = \D_2*_{\Z_2}\D_3$$

(cf. [MZ]). An interesting example of a maximal handlebody group which is not a
maximal bounded surface group is the second Mathieu group $M_{12}$.

\bigskip

1.4  {\it  Hurwitz actions versus maximal reducible actions}

\medskip

The smallest Hurwitz group is the linear fractional group 
${\rm PSL}_2(7)$  of order 168, acting on Klein's quartic ${\cal Q}_3$ of
genus 3. The smallest maximal handlebody group is the dihedral group $\D_6$
of genus 2, the second smallest the octahedral or symmetric group $\S_4$ of
genus 3 (see [Z5, Proposition 7]). The maximal handlebody group $\S_4$ is a
subgroup of index 7 in the Hurwitz group ${\rm PSL}_2(7)$, and its preimage
under the associated surjection
$\psi: (2,3,7) \to {\rm PSL}_2(7)$ is the
quadrangle group (2,2,2,3) (by an easy application of the formula of Riemann
Hurwitz). Since the restriction
$\psi:(2,2,2,3) \to \S_4$ factors through one of the three groups $G_2, G_3$
or $G_4$ (in fact through all of them), the action of $\S_4$
on Klein's quartic ${\cal Q}_3$ (unique up to conjugation) extends to a
handlebody $V_3$ of genus 3. Our first main result is the following.

\bigskip

{\bf Theorem 1.}   {\sl  Let $H$ be Hurwitz group of genus $g > 1$ acting
on a surface $\Sigma = \Sigma_g$.

\smallskip

i) If the action of $H$ on $\Sigma$ has a reducible 
subgroup $G$ of maximal possible order $12(g-1)$ then $\Sigma$ is a finite
regular covering of Klein's quartic and the
actions of $H$ and $G$ are obtained by lifting the actions of 
${\rm PSL}_2(7)$ and its subgroup $\S_4$ to $\Sigma$; in particular, $H$ surjects
onto ${\rm PSL}_2(7)$.  There are
infinitely many  Hurwitz actions with a maximal reducible subgroup.

\smallskip

ii) If $H$ is simple and not isomorphic to a subgroup of the alternating group
$\Bbb A_{24}$ (has no subgroup of index $\le 24$) then the action of every proper
subgroup of $H$ is reducible. }

\bigskip

In fact there seem to be very few cases of Hurwitz actions with a proper
irreducible subgroup,  see the Remark in section 4. Concerning part i) of Theorem
1, we note that  a lift of a subgroup $\S_4$ of ${\rm PSL}_2(7)$ to a finite regular
covering of Klein's quartic is always reducible but in general not completely
reducible;  at present we don't know if one can obtain infinitely many completely
reducible actions in this way.

\bigskip

{\bf Corollary 1.}   {\sl ${\rm PSL}_2(7)$ is the unique simple Hurwitz group which
realizes  both the smallest index reducible subgroup ($\Bbb S_4$ of index 7, with
preimage (2,2,2,3) in (2,3,7)) as well as
the largest  index irreducible subgroup ($\Bbb Z_7$ of index 24, with preimage
(7,7,7)).}

\bigskip

Considering just the types of groups and not the
actions, the following holds.

\bigskip

{\bf Theorem 2.}   {\sl A Hurwitz group which does not surject
onto ${\rm PSL}_2(7)$ is also a maximal reducible group; in particular, 
${\rm PSL}_2(7)$ is the only simple Hurwitz group which is not a maximal reducible
group. }

\medskip

The simple Hurwitz groups of order less than $10^6$ are of
linear fractional type ${\rm PSL}_2(q)$, and in addition the Janko group $J_1$
and the Hall-Janko group  $J_2$ (cf. [C]). 
The Hurwitz groups of linear fractional type are the following:

$q = 7$;

$q = p$ prime with $p \equiv \pm 1 {\; \rm mod \;} 7$ (each with a
unique Hurwitz action); 

$q = p^3$ with $p \equiv  \pm 2, \pm 3 {\; \rm mod \;} 7$ (each with three
different Hurwitz actions).

\bigskip

1.5 {\it Bounding and nonbounding Hurwitz actions}

\medskip

The action of the subgroup $\S_4$ of ${\rm PSL}_2(7)$ on Klein's quartic
${\cal Q}_3$ extends to a handlebody $V_3$ of genus 3; on the other hand, the
action of
${\rm PSL}_2(7)$ on  ${\cal Q}_3$ {\it does not bound}, i.e. does 
not extend to any compact 3-manifold
$M$ with exactly one boundary component $\partial M = {\cal Q}_3$ ([Z4,
Corollary 1]): the quotient orbifold  ${\cal Q}_3/{\rm PSL}_2(7)$ is
the 2-sphere with three branch points of orders 2, 3 and 7 which does not
occur as the unique boundary component of a compact 3-orbifold (since a
singular axis starting in the boundary point of order 7 can end only in a
dihedral point of type $\D_7$ but ${\rm PSL}_2(7)$ has no dihedral 
subgroup $\D_7$).

\medskip

The second smallest Hurwitz group is the linear fractional group
${\rm PSL}_2(8)$ of order 504 and genus 7, and it is shown in [Z4, Corollary
2] that the action by isometries of ${\rm PSL}_2(8)$ on the unique hyperbolic
surface  ${\cal Q}_7$ {\it bounds geometrically}, i.e. extends to a group of
isometries of a compact hyperbolic 3-manifold $M$ with totally geodesic
boundary $\partial M = {\cal Q}_7$. In particular, 
forgetting about the group action, the hyperbolic Hurwitz surface
${\cal Q}_7$ bounds geometrically. This leads naturally to the following:

\medskip

{\it Question} (see [Z6]):  Does Klein's quartic bound geometrically?

\medskip

Up to conjugation by homeomorphisms, there is a unique reducible action of
$\S_4$ on a surface of genus 3 (cf. [B]), represented by the subgroup $\Bbb S_4$ of
the Hurwitz action of ${\rm PSL}_2(7)$; we will show in section 3 that this
$\S_4$-action bounds geometrically, for some hyperbolic structure on a
surface of genus 3. But, in contrast to the irreducible Hurwitz actions, the
isometric action of $\S_4$ on a hyperbolic surface of genus 3 does not determine the
hyperbolic structure of the surface (since the deformation space of a
hyperbolic quadrangle of type (2,2,2,3)  has positive dimension whereas a
hyperbolic triangle of type (2,3,7) is unique up to isometry). So there are
uncountably many hyperbolic structures on a surface of genus 3 such that
$\S_4$ acts by isometries; only countably many of these can bound
geometrically, and at present we do not know if one can obtain Klein's
quartic by a construction as in section 3.

\medskip

So the reducible action of $\Bbb S_4$
extends both to a handlebody (in fact, $\Bbb S_4$ is the second smallest maximal
handlebody group), and also to a hyperbolic 3-manifold with totally geodesic
boundary; this raises the following:

\medskip

{\it Problem.} Show that every bounding action of a finite group on a closed
hyperbolic surface bounds also geometrically. In particular, if an action extends
to a handlebody, does it extend also to a hyperbolic 3-manifold with totally
geodesic boundary?

\medskip

For a proof of the following theorem, 
see [GZ, Corollary 3.14] (and [RZ, Proposition 2] for a proof that
${\rm PSL}_2(27)$ does not bound).

\bigskip

{\bf Theorem 3.}   {\sl  For $q \le 1000$, the only non-bounding Hurwitz actions of
linear fractional type ${\rm PSL}_2(q)$  are the Hurwitz actions of the groups
${\rm PSL}_2(7)$  and ${\rm PSL}_2(27)$; in fact all other Hurwitz actions 
with $q \le 1000$ bound
geometrically, and in particular the corresponding hyperbolic Hurwitz surfaces bound
geometrically.}

\bigskip

So there is some analogy between the Hurwitz actions of 
${\rm PSL}_2(7)$ and ${\rm PSL}_2(27)$: both do not bound, and in particular,
similar as for  ${\rm PSL}_2(7)$, it remains open whether the Hurwitz surface of the
action of ${\rm PSL}_2(27)$ bounds geometrically. Supported by some strong evidence
from  [GZ], we state the following:

\bigskip

{\it Conjecture.}  With the
exceptions of $q = 7$ and 27,  all Hurwitz actions of linear fractional type
${\rm PSL}_2(q)$  bound geometrically.

\bigskip

1.6  {\it  Maximal reducible groups versus maximal handlebody groups}

\medskip

Examples of small maximal handlebody groups which are not maximal bounded
surface groups are given in [Z3], and it is shown in [CZ] by extensive
computational methods  that 161 and 3781 are the two smallest genera for
which there exists a maximal handlebody group but not a maximal bounded surface
group (answering a question in [MZ]). In the next theorem, we confront maximal
reducible groups and maximal handlebody groups; the first claim is proved in [PZ,
Corollary 3.3].

\bigskip

{\bf Theorem 4.}   {\sl  i) The linear  fractional 
group ${\rm PSL}_2(q)$ is a maximal handlebody group exactly for all $q$ different
from 7, 9, 11 and  $3^{2m+1}$. The group ${\rm PSL}_2(27)$ is the smallest
simple group  which is  maximal reducible  but not a maximal handlebody
group. There are infinitely many maximal reducible groups which are not maximal
handlebody groups.

\smallskip

ii)  The smallest simple groups  which are maximal reducible, of order less than
979.200 (the order of the symplectic group ${\rm Sp}_4(4)$) and not of linear
fractional type, are the second Mathieu group
$M_{12}$, the Janko group $J_1$, the alternating group $\A_9$ and the 
Hall-Janko group $J_2$. The  group $M_{12}$ is a maximal handlebody group  of
types  $G_3$,  $G_4$ and $G_5$  but not a
maximal bounded surface group (not of type $G_2$),  
$J_1$ is a maximal bounded surface group, and 
$\Bbb A_9$  is a maximal   handlebody group  of type $G_3$  
but not a maximal bounded surface group. }

\medskip

So it remains open here whether ${\rm Sp}_4(4)$ is maximal reducible, and
whether the maximal reducible group  $J_2$ is a maximal handlebody or a
maximal  bounded surface group.

\bigskip

1.7 {\it  $G_7$-groups and bounding Hurwitz groups}

\medskip

An interesting class of maximal reducible groups closely related to Hurwitz
groups are the $G_7$-groups, 
$$G_7  =  \D_7*_{\Z_7}(2,3,7).$$
A surjection
$\varphi: G_7 \to G$ with torsionfree kernel 
 to a finite group $G$ determines a $G$-action  of maximal
possible order $12(g-1)$ on a product with handles $\cal P$ of genus $g$ (obtained
as the regular orbifold covering of the orbifold product with handles with
fundamental group $G_7$ associated to the kernel of
$\varphi$). Let $n$ be the index of the image $\varphi(2,3,7)$ of the triangle
group $(2,3,7)$ in $G$.

\medskip

If $n = 1$, the inner boundary $\partial_1 {\cal P}$ of  $\cal P$ consists of a single
surface of genus $g'$, with $84(g'-1) = 12(g-1)$, on which $G$ restricts to a 
Hurwitz action.  For arbitrary $n$, $\partial_1 {\cal P}$ consists of $n$ copies of
a Hurwitz surface, and the stabilizer in $G$ of each of these surfaces acts as a
Hurwitz group (isomorphic to $\varphi(2,3,7)$); see section 4 for an explicit
example.

\medskip

We say that a Hurwitz group {\it is bounding} if it admits a bounding Hurwitz action.

\bigskip

{\bf Theorem 5.}  {\sl  i) All bounding Hurwitz groups are $G_7$-groups. The groups
${\rm PSL}_2(7)$ and ${\rm PSL}_2(27)$  are the smallest non-bounding Hurwitz
groups, and  ${\rm PSL}_2(27)$  is the smallest simple non-bounding $G_7$-group.

\medskip

ii)  All Hurwitz  groups of linear fractional type except ${\rm PSL}_2(7)$  are $G_7$-groups. 
The  simple $G_7$-groups of
order less than $10^6$ are  Hurwitz groups of linear fractional type or one of   
the groups  ${\rm PSL}_2(49)$, $J_1$, $\A_9$ and $J_2$. The groups 
${\rm PSL}_2(49)$ and $\A_9$ are the smallest simple $G_7$-groups which are not
Hurwitz groups.}

\bigskip

Let $H$ be a Hurwitz group with a Hurwitz  action associated to a surjection
$\varphi: (2,3,7) \to H$. We denote by [2,3,7] the {\it extended triangle group}
generated by the reflections in the sides of a hyperbolic triangle with angles
$2\pi/2$, $2\pi/3$ and $2\pi/7$ (containing the triangle group (2,3,7) as a
subgroup of index 2,  of orientation-preserving elements). 
As a consequence of  [Z4, Theorem 3b] we have the following:

\bigskip

{\bf Theorem 6.}  {\sl  Let $H$ be a Hurwitz group with a Hurwitz action associated
to a surjection $(2,3,7) \to H$. If the surjection extends to a
surjection $[2,3,7] \to H$ then the Hurwitz action of $H$ bounds geometrically;
also, $H$ is a maximal bounded surface group.}

\bigskip

The Hurwitz groups $H$ as in Theorem 6 are exactly the "non-orientable Hurwitz
groups" (also called $H^*$-groups)  of maximal possible order
$84(g - 2)$  acting on a non-orientable surface of genus $g > 2$ (see [C2]). An
interesting example of such a group is the first Janko group
$J_1$, in particular a Hurwitz action of $J_1$ bounds geometrically (this
remains open for the Hall-Janko group $J_2$ which is also a Hurwitz group but not a
surjective image of the extended triangle group [2,3,7]). Since it is proved in [C3]
that, for $n \ge 168$, the alternating group
$\Bbb A_n$ is a surjective image of the extended triangle group [2,3,7], Theorem
6 implies:

\bigskip

{\bf Corollary 2.}  {\sl For each $n \ge 168$  there is a Hurwitz-action of the
alternating group $\Bbb A_n$ which bounds geometrically, in particular $\Bbb
A_n$ is a bounding Hurwitz group.}

\bigskip

The analogue remains open for the Hurwitz groups of linear fractional type ${\rm
PSL}_2(q)$ (but see the stronger Conjecture in section 1.5).

\bigskip

{\bf 2. Proofs of Theorems 1 and 2}

\medskip

In the notation of Theorem 1, the reducible group $G$ of order $12(g-1)$ has
index 7 in the Hurwitz group $H$ of order $84(g-1)$. There is a surjection 
$$\pi: (2,3,7) \to H$$
whose kernel is the universal covering group of the
surface $\Sigma_g$. Since the triangle group is a perfect group (has
trivial abelianization), also the Hurwitz group $H$ is perfect.  The preimage
$\pi^{-1}(G)$ of
$G$ in the triangle group (2,3,7) is a quadrangle group (2,2,2,3) (by an easy
application of the formula of Riemann-Hurwitz, this is the only signature of
a subgroup of index 7 in the triangle group (2,3,7)).

\medskip

The permutation
representation of $H$ by right multiplication on the seven
right cosets of $G$ in $H$ defines a nontrivial homomorphism  
$$\phi: H \to  \A_7$$
of  $G$ to the alternating group $\A_7$ of degree 7;
the image of $\phi$ is perfect and also a Hurwitz group (the 
kernel of the composition
$\phi \circ \pi: (2,3,7) \to \A_7$ is torsionfree since any normal subgroup
of (2,3,7) containing a torsion element is equal to (2,3,7)).
The perfect subgroups of the alternating group
$\A_7$ are the simple groups $\A_5$, $\A_6$, $\A_7$ and ${\rm PSL}_2(7)$ 
(see [CCN]). Since the only Hurwitz group among these groups is ${\rm
PSL}_2(7)$ (by [C] or by a direct computation using [GAP]), we have
surjections
$$\phi: H \to  {\rm PSL}_2(7),  \;\;\; \;\;\;
\phi \circ \pi: (2,3,7) \to {\rm PSL}_2(7),$$
and the second one gives the unique Hurwitz action
of ${\rm PSL}_2(7)$ on Klein's quartic ${\cal Q}_3$.

\medskip

The image $\phi(G)$ is a subgroup of index 1 or 7 in ${\rm
PSL}_2(7)$. Since ${\rm PSL}_2(7)$  is not a maximal handlebody group (not a
surjective image of the quadrangle group (2,2,2,3) as can easily be checked),
$\phi(G)$ has index 7 in ${\rm PSL}_2(7)$. Up to conjugation in 
${\rm PGL}_2(7)$, the only subgroup of index 7 in ${\rm PSL}_2(7)$ is the
symmetric group $\S_4$ of degree 4 (see [CCN]).

\medskip

The kernel  $U$  of the surjection  $\phi: H \to  {\rm PSL}_2(7)$ acts
freely on $\Sigma_g$ (since its preimage $\pi^{-1}(U)$ in the triangle group
(2,3,7) is torsionfree), and the actions of $H$ and $G$ on
$\Sigma_g$ project to the actions of ${\rm PSL}_2(7) \cong H/U$ and $\S_4
\cong G/U$ on Klein's quartic ${\cal Q}_3 = \Sigma_g/U$.

\medskip

In order to obtain infinitely many examples of such Hurwitz groups, one
considers finite-index characteristic subgroups of the fundamental group of
Klein's quartic and lifts the actions of ${\rm PSL}_2(7)$ and its subgroup
$\S_4$ to the corresponding finite regular coverings. 

\medskip

This concludes the proof of part i) of Theorem 1.  For the proof of part ii),
assume that $G$ is a proper subgroup of $H$, of index $n > 1$, which acts 
irreducibly on $\Sigma$.  Then the extension of $\pi_1(\Sigma)$ associated to
$G$ is a triangle group of index $n$ in the triangle group (2,3,7). The
triangle group of maximal possible index in (2,3,7) is the triangle group (7,7,7)
of index 24 which implies  $n \le 24$. The action by right multiplication of $G$ on
the right cosets of $H$ in $G$ defines a nontrivial, hence injective homomorphism
of the simple group  $H$ to the alternating group $\Bbb A_n$, in particular $H$
is isomorphic to a subgroup of $\Bbb A_{24}$.

\bigskip

For the {\it proof of Theorem 2}, suppose that the Hurwitz group $H$ does not
surject onto ${\rm PSL}_2(7)$. There is a surjection of the triangle group (2,3,7)
onto $H$, and the image of the quadrangle subgroup (2,2,2,3) of index 7 in (2,3,7)
has index 1 or 7 in $H$. If it has index 7 then, by the proof of Theorem 1, $H$
surjects onto ${\rm PSL}_2(7)$, otherwise the quadrangle group (2,2,2,3) surjects
onto $H$ and $H$ is a maximal reducible group.

\bigskip

Remark.  The Hurwitz action of ${\rm PSL}_2(7)$ has a cyclic irreducible subgroup
$\Bbb Z_7$; with this exception, every other cyclic subgroup of a
Hurwitz action of a simple group is reducible. In fact, the proper triangle
subgroups of the triangle group (2,3,7) are the groups (7,7,7), (3,7,7), (2,7,7)
and (3,3,7), of indices 24, 16, 9 and 8. The abelianizations of these
triangle groups are $\Bbb Z_7 \times \Bbb Z_7$, $\Bbb Z_7$, $\Bbb Z_7$ and $\Bbb
Z_3$, and only (7,7,7) admits a surjection with torsionfree kernel onto a cyclic
group ($\Bbb Z_7$); since (7,7,7) has index 24 in (2,3,7), the
maximal possible order of a simple  Hurwitz group with an irreducible cyclic
subgroup is $7 \cdot 24 = 168$, the order of the smallest Hurwitz group
${\rm PSL}_2(7)$. 

\medskip

The Hurwitz action of ${\rm PSL}_2(7)$ has also an irreducible  subgroup 
$\Bbb Z_7 \rtimes \Bbb Z_3$ of index 8, with preimage (3,3,7) in (2,3,7).
The second largest Hurwitz group ${\rm PSL}_2(8)$ has an irreducible subgroup 
$(\Bbb Z_2)^3 \rtimes \Bbb Z_7$ of index 9, with preimage (3,7,7) in (2,3,7).
We believe that these three examples are the only examples of
proper irreducible subgroups of a simple Hurwitz group (for example, the minimal
index of a proper subgroup of the third smallest Hurwitz group ${\rm PSL}_2(13)$ is
already 28, and hence the subgroup is reducible, and similar for all other Hurwitz
groups of linear fractional type). Here one should go over a list of the simple
(Hurwitz) groups of "small" order and consider the minimal index of a subgroup for
each group (see [C3] for the Hurwitz groups among the alternating groups).

\bigskip

Example.  The second largest Hurwitz group is the linear fractional
group   of order 504 and genus 7, and also the product 
${\rm PSL}_2(7) \times {\rm PSL}_2(8)$ is a Hurwitz group acting on a
unique surface of genus 1009 (by a direct computation using [GAP], see also
[C]). In fact, up to isomophism there is a unique surjection  
$\pi: (2,3,7) \to {\rm PSL}_2(7) \times {\rm PSL}_2(8)$. The
preimage  $\pi^{-1}(\S_4 \times {\rm PSL}_2(8))$ has index 7 in the triangle
group (2,3,7) and is a quadrangle group (2,2,2,3). We checked by [GAP]
that, for the corresponding surjection 
$\pi: (2,2,2,3) \to \S_4 \times {\rm PSL}_2(8)$, the image of the element
$x_1x_2 \in (2,2,2,3)$ has order $36$ and 36 is the minimal $n$ such
that the surjection factors through $G_n$. Hence the action of 
$\S_4 \times {\rm PSL}_2(8)$ on the surface of genus 1009 is reducible but
not completely reducible.  On the other hand, there are various surjections 
$(2,2,2,3) \to \S_4 \times {\rm PSL}_2(8)$ wich factor through $G_2$ or
$G_4$, so $\S_4 \times {\rm PSL}_2(8)$ is a maximal handlebody group.

\bigskip

{\bf 3.  Geometrically bounding actions}

\medskip

Up to isomorphisms there is a unique surjection of the
quadrangle group (2,2,2,3) onto $\S_4$; equivalently, up to conjugation by
homeomorphisms there is a unique action of $\S_4$ on a surface $\Sigma_3$ of
genus 3, with quotient orbifold  of type (2,2,2,3)  (cf. [B]).  There are
infinitely many hyperbolic structures on the quotient orbifold
$\Bbb H^2/(2,2,2,3)$; lifting such a hyperbolic structure to the surface
$\Sigma_3$, the group $\S_4$ acts by isometries.

\bigskip

{\bf Proposition 1.}   {\sl The unique $\S_4$-action on  a surface of
genus 3  with quotient orbifold of type 
(2,2,2,3)  bounds geometrically,
for some hyperbolic structure on the surface.}

\bigskip

{\it Proof.}  By Andreev's theorem ([T, chapter 13], [V, p.111])
there exists a unique finite hyperbolic polyhedron ${\cal C} = {\cal
C}(3,3,3)$  in hyperbolic 3-space
$\Bbb H^3$ which has the combinatorial structure of  a cube as follows. The
four edges of a face
$F$ of ${\cal C}$ have right dihedral angles (hence this face is orthogonal
to the four adjacent faces), the four edges connecting
$F$ to the opposite face $F'$ have three right dihedral angles and one dihedral
angle
$\pi/3$, and the face $F'$ has two consecutive edges with right dihedral
angles and two with dihedral angles $\pi/3$ such that the three edges of
${\cal C}$ with angles $\pi/3$ form a segment (note that the polyhedron
${\cal C}$ has no incompressible Euclidean 2-suborbifolds). Let $C$ denote
the orientation-preserving subgroup of index two in the group generated
by the reflections in the five faces of  ${\cal C}$ different form $F$ (so
$C$ is generated by rotations around the eight edges of ${\cal C}$ not in the
boundary of $F$); it is easy to check that there exists a
surjection from $C$ to $\S_4$ with torsionfree kernel.

\medskip

Let ${\cal D}$ denote the double of ${\cal C}$ along its face $F$ (again of
the combinatorial type of a cube, with five dihedral angles $\pi/3$), and 
by $D$ the orientation-preserving subgroup of index 2 in the Coxeter group
generated by the reflections in the faces of
${\cal D}$. There exists a surjection  $\pi: D \to \S_4$ with torsionfre
kernel $K$,   and hence $D/K \cong \S_4$ acts on the closed hyperbolic
3-manifold $M = \Bbb H^3/K$. The face $F$ of ${\cal C}$ lies in a hyperbolic
plane on which a quadrangle subgroup (2,2,2,3) of $D$ acts, and the
restriction of $\pi$ to this group (2,2,2,3) is also surjective.  The preimage
of the face $F$ of ${\cal C}$ is a closed totally geodesic surface 
$\Sigma_3$ of genus 3 in $M$ which is invariant under the action of
$\S_4$. The surface $\Sigma_3$ seperates $M$ into two isometric
hyperbolic 3-manifolds, each with an isometric action of $\S_4$, which have
the hyperbolic surface
$\Sigma_3$ as their common totally geodesic boundary, and hence the
isometric $\S_4$-action on the hyperbolic surface  $\Sigma_3$  bounds
geometrically, completing the proof of  Proposition 1.

\bigskip

We don't know if the hyperbolic surface 
$\Sigma_3$ is isometric to Klein's quartic (for example, one might try
to compute the lengths of the edges of the face $F$ of $Q$ and compare with
the edges of the corresponding quadrangle for a subgroup (2,2,2,3) of the
triangle group (2,3,7)).  Note that the construction can be modified by
choosing finite cubical polyhedra ${\cal C}$ with other small dihedral angles
for the face $F'$, obtaining other hyperbolic surfaces
$\Sigma_3$.

\medskip

There is a second action of $\Bbb S_4$ on a surface of genus 3 associated
to a surjection $(3,4,4)  \to  \Bbb S_4$ (see the list in [B] for the finite group
actions on a surface of genus 3), and this action is irreducible and bounds 
geometrically. In fact, it is a subgroup of an action of
$\Bbb S_4 \times \Bbb Z_2$ associated to a surjection $(2,4,6)  \to \Bbb S_4 \times
\Z_2$ which  bounds geometrically since this surjection extends to a surjection of
a hyperbolic tetrahedral group (the orientation-preserving subgroup of the group
generated by the reflections in the four faces of the tetrahedron)  associated to a
hyperbolic tetrahedron with one vertex at infinity of hyperbolic type (2,4,6)  
truncated by an othogonal plane, one edge of singular index 3 opposite to the edge
with singular index 6, and all other edges of singular index 2. The hyperbolic
surface of genus 3 determined by the action is not Klein's quartic since the
orientation-preserving isometry group 
${\rm PSL}_2(7)$ of Klein's quartic has no subgroup  $\Bbb S_4 \times \Z_2$. There
are exactly three other finite group actions on a surface of genus 3 of order
greater or equal to  48,  of orders 168, 96 and again 48, and each of these
three actions has a cyclic subgroup which does not bound, so we have the following:

\bigskip

{\bf Proposition 2.}   {\sl The unique bounding finite group action of largest
posible order on a surface of genus 3 is the action of $\Bbb S_4 \times \Z_2$
associated to a surjection $(2,4,6)  \to \Bbb S_4 \times \Z_2$.}

\bigskip

Since the reducible $\Bbb S_4$-action associated to a surjection 
$(2,2,2,3)  \to \Bbb S_4$ is also a subgroup of the action of $\Bbb S_4 \times
\Z_2$, this gives another proof of Proposition 1 (however, in contrast to
Proposition 1, we are sure here that this does not realize Klein's quartic).

\medskip

A description of the bounding and nonbounding finite group actions on a 
surface of genus 2 is given in [WZ]; for $g = 2$, no irreducible action bounds, and
the largest bounding action (extending to a handlebody) is by the smallest maximal
handlebody group $\Bbb D_6$.

\medskip

The situation for large group-actions on surfaces of genus 3 can be subsumed as
follows  (see [Z7] for some more details, and [B]  for a description of the groups
$D_{2,8,5}$ and  $D_{2,12,5}$).

\bigskip

{\bf Theorem 7.}   {\sl  The bounding and non-bounding finite group-actions on
a surface of genus 3, of order $\ge 24$, are the following.

\medskip
 
i) The two largest group-actions, 
of orders 168 and 96 are  represented by surjections
$$(2,3,7) \to  {\rm PSL}_2(7) \;\; {\sl and} \;\; (2,3,8)  \to  
\D_3 \ltimes (\Z_4 \times \Z_4)$$

and do not bound. 

\medskip

ii) Two actions of order 48 associated to surjections
$$(3,3,4)  \to \Z_3 \ltimes (\Z_4 \times \Z_4) \;\; {\it and} \;\;
(2,4,6) \to  \Z_2 \times S_4;$$ 

the first one is a subgroup of index 2 of the group of order 96 in i) and does not
bound, the second one is the largest bounding group-action on a surface of genus 3;
it bounds geometrically but does not extend to  a handlebody.

\medskip

iii) Two non-bounding actions of order 32 associated to surjections
$$(2,4,8)  \to \Z_2 \ltimes (\Z_2 \times \Z_8) \;\; {\it and} \;\;
(2,4,8)  \to \Z_2 \times D_{2,8,5}.$$

\medskip

iv) Two non-bounding actions of order 24  associated to surjections

$$(3,3,6) \to  {\rm SL}_2(3)  \;\; {\it and} \;\;  (2,4,12) \to  D_{2,12,5}.$$

\medskip

v)  Three bounding actions of order 24 associated to surjections

$$(2,6,6) \to \Z_2 \times \A_4, \;\;  (3,4,4) \to \S_4  \;\; {\it and} \;\; 
(2,2,2,3) \to \S_4;$$ 

these three actions are 
subgroups of index 2 of the geometrically bounding action of order 48 in ii). The
last one is also the largest action on a surface of genus 3 which extends to a
handlebody (in fact $\S_4$ is the unique maximal handlebody group of genus 3,  
of maximal possible order $12(g-1)$).

\medskip

Finally, for each of the non-bounding actions in i) - iv) there is already a
cyclic subgroup which does not bound. }

\bigskip

A beautiful example of a geometrically bounding action
on a surface of genus 4 is the $\Bbb A_5$-action described in [Z6] associated to a
surjection $(2,5,5) \to \Bbb A_5$. This is a subgroup of an action of $\Bbb S_5$
associated to a surjection $(2,4,5)  \to  \Bbb S_5$ which realizes the maximal order
of a finite group action on a surface of genus 4 and is also 
geometrically bounding.

\bigskip

{\bf Proposition 3.}   {\sl The largest group-action on a surface of genus 4, of
type 
$(2,4,5)  \to \S_5$, bounds geometrically. The largest action in genus 4 which
extends to a handlebody is of type   $(2,2,2,3)  \to \D_3 \times \D_3$.}

\bigskip

As noted
before, the first (geometrically) bounding Hurwitz actions is that of ${\rm
PSL}_2(8)$ on a surface of genus 7.

\bigskip

{\bf 4. Proofs of Theorems 4 and 5}

\medskip

As noted before, the first claim of Theorem 4 is proved in  [PZ, Corollary 3.3] 
(note that the case $q = 11$ was overlooked in [PZ].) In particular,  ${\rm
PSL}_2(27)$ is not a maximal handlebody group (which can be checked easily also by
a direct computation). The unique Hurwitz action of 
${\rm PSL}_2(27)$ is associated to a surjection $(2,3,7) \to {\rm PSL}_2(27)$; 
since ${\rm PSL}_2(27)$ has subgroups $\D_7$, this extends to a surjection 
$G_7  =  \D_7*_{\Z_7}(2,3,7)  \to {\rm PSL}_2(27)$, hence ${\rm PSL}_2(27)$
is a $G_7$-group.

\medskip

The group $G = {\rm PSL}_2(27)$ is a Hurwitz group in a unique way: up to
equivalence there is a unique surjections $(2,3,7) \to G$; let
$K$ denote the kernel of this surjection which is a surface group of genus
$g = 118$. By abelianization and reduction of coordinates mod $p$, for some
prime $p$, we get a surjection $K \to (\Bbb Z_p)^{2g}$ whose kernel is a
charateristic subgroup $\tilde K$ of $K$, in particular normal in $(2,3,7)$,
so we have an exact sequence

$$ 1 \to \tilde K  \hookrightarrow (2,3,7) \to (\Bbb Z_p)^{2g}
\rtimes {\rm PSL}_2(27)  \to 1.$$

Since $G = {\rm PSL}_2(27)$ has  dihedral subgroups $\Bbb D_7$ (of order 14), the
surjection
$(2,3,7)  \to (\Bbb Z_p)^{2g} \rtimes G$ extends to a surjection with
torsionfree kernel of
$G_7 =  \D_7*_{\Z_7}(2,3,7)$ to $(\Bbb Z_p)^{2g} \rtimes G$, hence 
$(\Bbb Z_p)^{2g} \rtimes G$ is a maximal reducible $G_7$-group of Hurwitz
type. On the other hand,  $(\Bbb Z_p)^{2g} \rtimes G$ is not a maximal
handlebody group for $p > |G|$ since a surjection with
torsionfree kernel $G_n \to  (\Bbb Z_p)^{2g} \rtimes G$ would induce a surjection
$G_n \to G$ with torsionfree kernel which does not exist for $2 \le n \le 5$ (since
$G = {\rm PSL}_2(27)$ is not a maximal handlebody group).

\medskip

This completes the proof of Theorem 4 i). 
Theorem 4 ii) and Theorem 5 ii) can easily be
checked by case-by-case computations (using e.g. [GAP]), going over the list of the
simple groups of small order, their maximal subgroups and the orders and conjugacy
classes of elements (see [CNN]). The first statement of Theorem 5 ii) follows
from the well-known classification of the subgroups of the linear fractional groups 
${\rm PSL}_2(q)$; in particular, with the exception of ${\rm PSL}_2(7)$, an
element of order 7 in a linear fractional Hurwitz group $H$ lies in a dihedral
subgroup $\Bbb D_7$ of order 14, and hence $H$ is a $G_7$-group. Theorem 5 i) follows
from the argument given in section 1.5 for the case of ${\rm PSL}_2(7)$.

\bigskip

Example.  The Hurwitz group ${\rm PSL}_2(7)$ of genus 3 is
not a $G_7$-group since an element of order 7 does not
lie in a dihedral subgroup $\Bbb D_7$.  We finish with an example of a maximal
reducible  action of ${\rm PGL}_2(7)$ on a surface
$\Sigma$ of genus 29 and the induced action on a product with handles
$\cal P$ with outer boundary $\partial_0 {\cal P} = \Sigma$ whose inner boundary 
$\partial_1 {\cal P}$  consists of two copies of Klein's quartic with the Hurwitz
action of ${\rm PSL}_2(7)$.

\medskip

The kernel of the composition of surjections

$$(2,2,2,3)  \to  \D_7*_{\Z_7}(2,3,7) \to {\rm PGL}_2(7)$$

is a surface group $\pi_1(\Sigma)$ of genus $g = 29$ (note that, in contrast to 
${\rm PSL}_2(7)$,  ${\rm PGL}_2(7)$ has dihedral subgroups $\D_7$). 
Let $\bar {\cal P}$ denote the orbifold product with handles with orbifold
fundamental group  $\D_7*_{\Z_7}(2,3,7)$ (cf. [R]), and let $\cal P$ be the product
with handles of genus $g = 29$ 
which is the covering of $\bar {\cal P}$ associated to the kernel of the surjection
$\D_7*_{\Z_7}(2,3,7) \to {\rm PGL}_2(7)$. The outer boundary $\partial_0 {\cal P}$
is the surface  $\Sigma$; since 
the triangle group $(2,3,7)$ maps
onto the subgroup ${\rm PSL}_2(7)$ of index 2 in ${\rm PGL}_2(7)$, the inner
boundary 
$\partial_1 {\cal P}$ consists of two copies of Klein's quartic of genus 3
with the Hurwitz action of ${\rm PSL}_2(7)$.

\bigskip  \bigskip

\centerline {\bf References}

\medskip

\item {[B]} S.A. Broughton, {\it Classifuing finite group actions on
surfaces of low genus.}  J. Pure Appl. Alg. 69 (1990),  233-270

\item {[C1]} M. Conder, {\it The genus of compact Riemann surfaces with
maximal automorphism group.}  J. Algebra 108 (1987),  204-247

\item {[C2]} M. Conder, {\it Maximal automorphism groups of symmetric 
Riemann surfaces with small genus.}  J. Algebra
114 (1988), 16-28

\item {[C3]} M. Conder, {\it Generators for alternating and symmetric groups.} 
J. London Math. Soc. 22  (1980),  75-86

\item {[CZ]} M.D.E. Conder, B. Zimmermann,  {\it Maximal bordered surface
groups versus maximal handlebody groups.}  Contemporary Math.  629 (2014),
99-105

\item {[CCN]} J.H. Conway, R.T. Curtis, S.P. Norton, R.A. Parker, R.A.
Wilson,  {\it Atlas of Finite Groups.} Oxford University Press 1985

\item {[GAP]}  The GAP Group, GAP - Groups, Algorithms, and Programming, 
Computational algebra system available under http://www.gap-system.org

\item {[GZ]} M. Gradolato, B. Zimmermann,  {\it  Extending finite group
actions on  surfaces to hyperbolic 3-manifolds.} Math. Proc. Cambridge
Phil. Soc.  117   (1995), 137-151

\item {[MMZ]} D. McCullough, A. Miller, B. Zimmermann,  {\it Group actions on
handlebodies.}  Proc. Lond. Math. Soc. 59 (1989), 373-415

\item {[MY]} W.H. Meeks, S.T. Yau,  {\it The equivariant Dehn's lemma and
loop theorem.}  Comment. Math. Helv. 56 (1981), 225-239

\item {[MZ]}  A. Miller, B. Zimmermann,  {\it  Large groups of symmetries of
handlebodies.} Proc. Amer. Math. Soc. 106  (1989), 829-838

\item {[PZ]}  L. Paoluzzi, B. Zimmermann  {\it  On finite quotients of the
Picard group and related hyperbolic tetrahedral groups.}
Rend. Istit. Mat. Univ. Trieste Suppl. 1 Vol. 32, 257-288 (2001)

\item {[R]}  M. Reni,  {\it  A graph-theoretical approach to Kleinian groups.} 
Proc. London Math. Soc. 67  (1993), 200-224

\item {[RZ]}  M. Reni, B. Zimmermann,  {\it  Extending finite group actions from
surfaces to handlebodies.}  Proc. Amer. Math. Soc. 124  (1996), 2877-2887

\item {[T]} W. Thurston, {\it The geometry and topology o 3-manifolds.} 
Lecture Notes, Dept. of Math., Princeton Univ., Princeton 1977

\item {[V]} E.B. Vinberg (Ed.), {\it Geometry II.}  Encyclopaedia of Math.
Sciences vol. 29, Springer-Verlag 1993

\item {[WZ]} C. Wang, S. Wang, Y. Zhang, B. Zimmermann,  {\it Finite group
actions on the genus-2 surface, geometric generators and extendability.}  
Rend. Istit. Mat. Univ. Trieste 52 (2020), 513-524

\item {[Z1]} B. Zimmermann,  {\it \"Uber Abbildungsklassen von
Henkelk\"orpern.}  Arch. Math. 33 (1979), 379-382

\item {[Z2]} B. Zimmermann,  {\it \"Uber Hom\"oomorphismen n-dimensionaler
Henkelk\"orper und endliche Erweiterungen von Schottky-Gruppen.}
Comm. Math. Helv. 56 (1981), 474-486

\item {[Z3]}  B. Zimmermann,  {\it Finite group actions on handlebodies and
equivariant Heegaard genus for 3-manifolds.} Topol. Appl. 43, 263-274 (1992)

\item {[Z4]}  B. Zimmermann, {\it  Hurwitz groups and finite group actions on
hyperbolic 3-manifolds.}  J. London Math. Soc. 52  (1995), 199-208

\item {[Z5]} B. Zimmermann,  {\it Genus actions of finite groups on
3-manifolds.}  Michigan Math. J. 43 (1996),  593-610

\item {[Z6]}  B. Zimmermann, {\it  A note on surfaces bouonding hyperbolic
3-manifolds.}  Monatsh. Math. 142  (2004), 267-273

\item {[Z7]}  B. Zimmermann, {\it  A note on large bounding and non-bounding
finite group-actions on surfaces of small genus.}  arXiv:2205.14425

\bye